\documentclass{amsart}
\usepackage{amsmath,amssymb}
\usepackage{graphicx}
\usepackage{xcolor}
\usepackage{hyperref}

\newtheorem{theorem}{Theorem}

\newtheorem{remark}{Remark}

\begin{document}

\title[Scale-Local Deformations of Vortex Rings in Euler Flows]
{Nonlinear Scale-Local Geometric Deformations of Vortex Rings in Smooth Euler Flows via Bayesian Optimization and Adjoint Methods} 



\author{Tsuyoshi Yoneda} 
\address{Graduate School of Economics, Hitotsubashi University, 2-1 Naka, Kunitachi, Tokyo 186-8601, Japan} 
\email{t.yoneda@r.hit-u.ac.jp}


\date{\today} 


\begin{abstract} 
We consider the incompressible three-dimensional Euler equations for a vortex ring with Kelvin waves undergoing radially expanding Lagrangian transport.
To clarify the fundamental mechanisms underlying nonlinear scale-local deformations of the vortex structure, we develop a geometric Lagrangian framework that avoids singular integral representations of the pressure and yields a novel wave equation governing the axis of swirling particles.
Within this framework, we identify intrinsic nonlinear mechanisms that drive scale-local deformations of the vortex structure, supported by a machine-learning-based analysis.
Specifically, we propose a hybrid optimization framework that combines Bayesian global exploration with adjoint-based local refinement.
The resulting optimization problem exhibits a highly non-convex loss landscape, in which the adjoint method alone fails to escape local minima.
\end{abstract}

\maketitle

\section{Introduction}

One of the fundamental features of Navier–Stokes turbulence is that it is not purely random, but instead organized by vortex-stretching dynamics.
More precisely, direct numerical simulations at sufficiently high Reynolds numbers \cite{Goto-2008,GSK,Motoori-2019,Motoori-2021} have revealed the existence of a hierarchical structure of vortex-stretching motions.
In particular, Goto–Saito–Kawahara \cite{GSK} demonstrated that turbulence in a periodic cube exhibits a self-similar hierarchy of antiparallel vortex-tube pairs, sustained by the generation of smaller-scale vortices through stretching within larger-scale strain fields.
This observation was further explored by Y–Goto–Tsuruhashi \cite{YGT} and Tsuruhashi–Goto–Oka–Y \cite{TGOY} (see also \cite{JNY, JY1,JY2,JY3} for related mathematical results).
Following this line of research, attention naturally turns to scale-local deformations of vortical structures, which play a central role in the energy cascade.

A natural setting in which such mechanisms can be examined in a controlled manner is provided by vortex-ring collisions.
Many studies have investigated this problem (see \cite{MOMPBR} and references therein), and the main findings can be summarized as follows.
When two vortex rings of equal strength and opposite circulation collide head-on, a strong strain field is generated.
The dominant instability mechanism depends on the Reynolds number.
At low Reynolds numbers, long-wavelength Crow-type instabilities prevail, leading to relatively mild deformations.
In contrast, at high Reynolds numbers, short-wavelength perturbations at the scale of the vortex core radius become dominant, and the dynamics are governed by the elliptical instability.
However, since the elliptical instability is derived from a linear stability analysis, it does not by itself provide a rigorous description of the nonlinear evolution.

Motivated by these considerations, we develop a geometric Lagrangian framework for the three-dimensional incompressible Euler equations that avoids singular integral representations of the pressure.
Within this framework, we rigorously derive a novel wave equation governing the evolution of the axis of swirling particles (denoted by $\zeta$, see \eqref{zeta}), and employ a machine-learning-based analysis to identify intrinsic nonlinear mechanisms responsible for the emergence of scale-local deformations of the vortex structure (the
corresponding Python code is available in the
repository \cite{Y}).
More precisely, when the axis of swirling particles and the vortex axis (denoted by $e_s$, see \eqref{es}) are not aligned, the resulting rotating flow tends to exhibit shear.
From a physical viewpoint, it is therefore natural to expect that the system favors configurations in which these axes remain aligned during the evolution, thereby minimizing shear.
Motivated by this observation, we introduce the mean square directional correlation (MSDC) as a functional (see \eqref{MSDC}), which we adopt as an objective functional to be maximized.

\section{Derivation of Novel Wave Equations}

In this section we derive novel wave equations
that exhibit swirling particles near the vortex axis.
Let $\Phi$ denote the Lagrangian flow map associated with the incompressible Euler equations, written in Lagrangian form (omitting the explicit volume-preserving constraint) as follows:
\begin{equation}\label{Euler eq.}
\left(\partial_t^2 \Phi(t,\cdot)\right)\circ \Phi^{-1}(t,x)
= - \nabla p(t,x),
\qquad t>0,\ x\in\mathbb{R}^3.
\end{equation}
Let $\ell=\bigcup_s \ell(s)\subset\mathbb{R}^3$ be a smooth curve representing the vortex axis.
We now describe the time evolution of swirling particles around $\ell$.
For fixed $s$, let $z$ denote the arc-length parameter along the transported curve
$\bigcup_t \Phi(t,\ell(s))$.
More precisely, we define a reparametrization $t=t(z)$ (implicitly depending on $s$) by
\begin{equation}\label{ODE}
\partial_z t = v(t)^{-1}
\qquad
v(t):=\left|\partial_t \Phi(t,\ell(s))\right|,
\end{equation}
so that
\[
\left|\partial_z \Phi(t(z),\ell(s))\right|=1.
\]
Note that this implies $\partial_z = v(t)^{-1}\partial_t$.
Using the transported curve $\bigcup_t \Phi(t,\ell(s))$, we introduce a Lagrangian moving frame.
First, we define the unit tangent vector
\[
\tau(z):=\partial_z \Phi(t(z),\ell(s)).
\]
The curvature $\kappa(z)$ and unit normal vector $n(z)$ are defined by
\[
\partial_z \tau(z)=\kappa(z)\,n(z).
\]
The torsion $T(z)$ and binormal vector $b(z)$ are then defined through the Frenet--Serret relation:
\[
\partial_z n(z)=-\kappa(z)\,\tau(z)+T(z)\,b(z).
\]
For any initial particle
\[
x=\ell(s)
+Z(0,r_1,r_2)\,\tau(0)
+r_1\,n(0)+r_2\,b(0),
\]
the flow map $\Phi(t,x)$ can be uniquely expressed (for sufficiently small $r_1,r_2$) as
\begin{equation}\label{Phi}
\begin{split}
\Phi(t,x)
&=:\Phi(t,s,r_1,r_2)\\
&=\Phi(t,\ell(s))
+Z(t,r_1,r_2)\,\tau(t)
+R_1(t,r_1,r_2)\,n(t)
+R_2(t,r_1,r_2)\,b(t).
\end{split}
\end{equation}
Here
\[
Z(0,r_1,r_2)=\alpha_1(0)r_1+\alpha_2(0)r_2,
\qquad
Z(t,0,0)=0
\]
and
\[
R_1(0,r_1,r_2)=r_1,
\qquad
R_2(0,r_1,r_2)=r_2
\]
with $\alpha_1(0),\alpha_2(0)\in\mathbb{R}$.
Since the Jacobian
\[
\frac{\partial(R_1,R_2)}{\partial(r_1,r_2)}
\]
is nonzero at $(r_1,r_2)=(0,0)$, the inverse function theorem implies that, for each fixed $t>0$ and for sufficiently small $|R_1|,|R_2|$, we may rewrite \eqref{Phi} in terms of $(R_1,R_2)$ as follows.
Let $r_1=r_1(t,R_1,R_2)$ and $r_2=r_2(t,R_1,R_2)$ denote the corresponding inverse functions.
Then, for any initial particle
\[
x=\ell(s)
+r_1(t,R_1,R_2)\,n(0)
+r_2(t,R_1,R_2)\,b(0),
\]
we have
\begin{equation}\label{explicit-formula}
\begin{split}
\Phi(t,x)
&=:\Phi(t,s,R_1,R_2)\\
&=\Phi(t,\ell(s))
+Z\bigl(t,r_1(t,R_1,R_2),r_2(t,R_1,R_2)\bigr)\,\tau(t)\\
&\quad +R_1\,n(t)+R_2\,b(t).
\end{split}
\end{equation}
Note that $\Phi(t,s,R_1,R_2)\big|_{R_1=R_2=0}$ represents the Lagrangian transport of the vortex axis $\ell$.
We now simplify the function $Z$ by the expansion
\begin{equation}\label{Z}
Z\bigl(t,r_1(t,R_1,R_2),r_2(t,R_1,R_2)\bigr)
=
\alpha_1(t)R_1+\alpha_2(t)R_2+\mathcal{O}(R_1^2+R_2^2),
\end{equation}
and we interpret $Z$ as describing the time evolution of swirling particles around the vortex axis.
We are now in a position to state the following key theorem.
\begin{theorem}\label{criterion}
The quantities $\alpha_1$ and $\alpha_2$
satisfy the following second--order wave equations:
\begin{equation}\label{wave equations}
\begin{cases}
\displaystyle
\alpha_1''(t)
=\frac{v''(t)}{v(t)}\,\alpha_1(t)
+2v(t)\kappa'(t)
+4v'(t)\kappa(t),\\
\displaystyle
\alpha_2''(t)
=\frac{v''(t)}{v(t)}\,\alpha_2(t).
\end{cases}
\end{equation}
\end{theorem}

\begin{remark}
In a separate line of research, Hasimoto's remarkable discovery that the
Local Induction Approximation (LIA) reduces to the focusing nonlinear
Schr\"odinger equation revealed a profound connection between filament
geometry and dispersive wave dynamics (see \cite{H}).
Despite these influential developments, existing theories rely fundamentally
on vortex--filament approximations, such as concentrated vorticity,
Biot--Savart expansions, LIA truncations, or asymptotic regimes in which the
vortex core radius is assumed to be small.
Consequently, the emergence of wave phenomena has traditionally been viewed
as a feature of singular or nearly singular vortex models.
It has therefore remained unclear whether similar wave dynamics arise
intrinsically within the full incompressible Euler equations in a smooth,
non--filament setting.

To the best of our knowledge, this paper is the first derivation showing that wave 
equations governing Lagrangian motion emerge directly from the 
incompressible Euler equations in a fully smooth setting. 
\end{remark}

\section{Proof of Theorem \ref{criterion}}\label{proof}
The key point of the proof is a direct computation of
$\partial_{R_1}(\partial_t^2\Phi\cdot\tau)$ and $\partial_{R_2}(\partial_t^2\Phi\cdot\tau)$
along the transported curve $\bigcup_t \Phi(t,\ell(s))$.
Let us recall the Frenet--Serret formulas:
\begin{equation*}
\frac{d}{dz}
\begin{pmatrix}
\tau\\
n\\
b
\end{pmatrix}
=
\begin{pmatrix}
0&\kappa&0\\
-\kappa& 0& T\\
0& -T& 0
\end{pmatrix}
\begin{pmatrix}
\tau\\
n\\
b
\end{pmatrix}.
\end{equation*}
By combining the Frenet--Serret formulas with \eqref{explicit-formula} and \eqref{Z}, we obtain the expansions
\begin{equation}\label{Frenet-Serret}
\begin{cases}
\partial_{z}\Phi
=
\tau+R_1(Tb-\kappa \tau)-R_2 T n\\
\qquad\quad +(R_1\alpha_1+R_2\alpha_2)\kappa n+(R_1\partial_z\alpha_1+R_2\partial_z\alpha_2)\tau
+\mathcal O(R_1^2+R_2^2),\\[2mm]
\partial_{R_1}\Phi
=
n+\alpha_1\tau+\mathcal O(|R_1|+|R_2|),\\[1mm]
\partial_{R_2}\Phi
=
b+\alpha_2\tau+\mathcal O(|R_1|+|R_2|),
\end{cases}
\end{equation}
where $\mathcal O(\cdot)$ denotes Landau's notation.
Equivalently, \eqref{Frenet-Serret} can be written in matrix form as
\begin{equation*}
\begin{split}
\begin{pmatrix}
\partial_{z}\Phi\\
\partial_{R_1}\Phi\\
\partial_{R_2}\Phi
\end{pmatrix}
&=
\begin{pmatrix}
(1-\kappa R_1)+(R_1\partial_z\alpha_1+R_2\partial_z\alpha_2)
&
-R_2T +(R_1\alpha_1+R_2\alpha_2)\kappa
&
R_1 T\\
\alpha_1& 1& 0\\
\alpha_2& 0& 1
\end{pmatrix}
\begin{pmatrix}
\tau\\
n\\
b
\end{pmatrix}\\
&=:\
\begin{pmatrix}
A&B&C\\
\alpha_1& 1& 0\\
\alpha_2& 0& 1
\end{pmatrix}
\begin{pmatrix}
\tau\\
n\\
b
\end{pmatrix},
\end{split}
\end{equation*}
where
\begin{equation*}
\begin{cases}
A=(1-\kappa R_1)+(R_1\partial_z\alpha_1+R_2\partial_z\alpha_2),\\
B= -R_2T +(R_1\alpha_1+R_2\alpha_2)\kappa,\\
C= R_1 T.
\end{cases}
\end{equation*}
We omit higher--order terms since we will eventually take the limit $R_1,R_2\to0$.
The corresponding inverse matrix is
\begin{equation}\label{inverse matrix}
\begin{pmatrix}
\tau\\
n\\
b
\end{pmatrix}
=
\frac{1}{D}
\begin{pmatrix}
1& -B& -C\\
-\alpha_1& A-C\alpha_2& C\alpha_1\\
-\alpha_2& B\alpha_2& A-B\alpha_1
\end{pmatrix}
\begin{pmatrix}
\partial_{z}\Phi\\
\partial_{R_1}\Phi\\
\partial_{R_2}\Phi
\end{pmatrix},
\end{equation}
where $D=A-B\alpha_1-C\alpha_2$.
Moreover, differentiating \eqref{Frenet-Serret} once more in $z$ and using the Frenet--Serret relations, we obtain
\begin{equation*}
\begin{split}
\partial_z^2\Phi
=&\ \kappa n
-R_1(T^2+\kappa^2 )n
+R_1\bigl((\partial_zT)b-(\partial_z\kappa) \tau\bigr)\\
&\ -R_2T(-\kappa\tau+Tb)
-R_2(\partial_zT) n\\
&\ +(R_1\partial_z^2\alpha_1+R_2\partial_z^2\alpha_2)\tau
+2(R_1\partial_z\alpha_1+R_2\partial_z\alpha_2)\kappa n\\
&\ +(R_1\alpha_1+R_2\alpha_2)(\partial_z\kappa) n
+(R_1\alpha_1+R_2\alpha_2)\kappa(-\kappa\tau+Tb)\\
&\ +\mathcal O(R_1^2+R_2^2).
\end{split}
\end{equation*}
Consequently,
\begin{equation*}
\begin{split}
\partial_{R_1}\partial_z\Phi\big|_{R_1=R_2=0}
&=Tb-\kappa \tau+\alpha_1\kappa n+\partial_z\alpha_1\,\tau,\\
\partial_{R_1}\partial_z^2\Phi\big|_{R_1=R_2=0}
&=-(T^2+\kappa^2)n+\bigl((\partial_zT)b-(\partial_z\kappa)\tau\bigr)\\
&\qquad
+\partial_z^2\alpha_1\,\tau+2(\partial_z\alpha_1)\kappa n
+\alpha_1(\partial_z\kappa) n-\alpha_1\kappa(\kappa\tau-Tb),\\
\partial_{R_2}\partial_z\Phi\big|_{R_1=R_2=0}
&=-Tn+\alpha_2\kappa n+\partial_z\alpha_2\,\tau,\\
\partial_{R_2}\partial_z^2\Phi\big|_{R_1=R_2=0}
&=-T(-\kappa \tau+Tb)-(\partial_zT)n\\
&\qquad+\partial_z^2\alpha_2\,\tau+2(\partial_z\alpha_2)\kappa n
+\alpha_2(\partial_z\kappa)n-\alpha_2\kappa(\kappa\tau-Tb).
\end{split}
\end{equation*}
Thus, using the identities obtained above and the orthonormality of the Frenet frame, we have
\begin{equation*}
\bigl(\partial_{R_1}\partial_z^2\Phi\bigr)\cdot\tau\Big|_{R_1=R_2=0}
=-\partial_z\kappa+\partial_z^2\alpha_1-\kappa^2\alpha_1,
\qquad
\bigl(\partial_{R_1}\partial_z\Phi\bigr)\cdot\tau\Big|_{R_1=R_2=0}
=-\kappa+\partial_z\alpha_1,
\end{equation*}
and
\begin{equation*}
\bigl(\partial_{R_2}\partial_z^2\Phi\bigr)\cdot\tau\Big|_{R_1=R_2=0}
=T\kappa+\partial_z^2\alpha_2-\kappa^2\alpha_2,
\qquad
\bigl(\partial_{R_2}\partial_z\Phi\bigr)\cdot\tau\Big|_{R_1=R_2=0}
=\partial_z\alpha_2.
\end{equation*}
On the other hand, since $\partial_t = v\,\partial_z$ along the vortex axis, the Leibniz rule gives
\begin{equation*}
\partial_t^2\Phi
=
v^2\partial_z^2\Phi+v'\partial_z\Phi,
\end{equation*}
where $v'=\partial_t v$ (and we omit writing the variables for brevity). Recall that the prime notation ${}'$ denotes differentiation with respect to time $t$, and that $\partial_z$ denotes differentiation with respect to the arc--length parameter $z$.
Therefore, for $R_1=R_2=0$, we obtain
\begin{equation}\label{pressure-estimate}
\begin{split}
-\partial_{R_1}(\nabla p\cdot\tau)
&=
\partial_{R_1}(\partial_t^2\Phi\cdot\tau)
=(\partial_{R_1}\partial_t^2\Phi)\cdot\tau\\
&=
(-\kappa+\partial_z\alpha_1)v'
+(-\partial_z\kappa+\partial_z^2\alpha_1-\kappa^2\alpha_1) v^2\\
&=
-v'\kappa -v^2\partial_z\kappa+\partial_t^2\alpha_1-v^2\kappa^2\alpha_1,\\[2mm]
-\partial_{R_2}(\nabla p\cdot\tau)
&=
\partial_{R_2}(\partial_t^2\Phi\cdot\tau)
=(\partial_{R_2}\partial_t^2\Phi)\cdot\tau\\
&=
(\partial_z\alpha_2)v'
+(T\kappa+\partial_z^2\alpha_2-\kappa^2\alpha_2)v^2\\
&=
v^2T\kappa+\partial_t^2\alpha_2-v^2\kappa^2\alpha_2.
\end{split}
\end{equation}
Next, we derive complementary identities using the Euler equations.
Along the vortex axis we have $\partial_t\Phi = v\tau$, and hence
\begin{equation*}
\kappa n=\partial_z^2\Phi(t(z),\ell(s))
=\partial_z\bigl(\partial_t\Phi\, v^{-1}\bigr)
=\partial_t^2\Phi\, v^{-2}-\partial_t\Phi\, v'\, v^{-3}.
\end{equation*}
Consequently,
\begin{equation*}
\partial_t^2\Phi=v^2\kappa n+v'\tau.
\end{equation*}
Using the Euler equations $\partial_t^2\Phi=-(\nabla p)\circ\Phi$, we obtain, for $R_1=R_2=0$ (again omitting $\circ\Phi$),
\begin{eqnarray*}
-\nabla p\cdot\tau &=& \partial_t^2\Phi\cdot\tau = v',\\
-\nabla p\cdot n    &=& \partial_t^2\Phi\cdot n    = v^2\kappa,\\
-\partial_z(\nabla p\cdot n) &=& \partial_z(v^2\kappa)=v^2\partial_z\kappa+2v'\kappa,\\
-\nabla p\cdot b &=& \partial_t^2\Phi\cdot b = 0.
\end{eqnarray*}
By the explicit inverse formula \eqref{inverse matrix} and the expansion
\begin{equation*}
\begin{split}
D^{-1}
=
1&-\Bigl((-\kappa+\partial_z\alpha_1)-\alpha_1^2\kappa-T\alpha_2\Bigr)R_1\\
&-\Bigl(\partial_z\alpha_2-(-T+\alpha_2\kappa)\alpha_1\Bigr)R_2
+\mathcal O(R_1^2+R_2^2),
\end{split}
\end{equation*}
we can write
\begin{equation*}
\begin{split}
\nabla p\cdot \tau
=&\
D^{-1}(\nabla p\cdot \partial_z\Phi)
-BD^{-1}(\nabla p\cdot\partial_{R_1}\Phi)
-CD^{-1}(\nabla p\cdot \partial_{R_2}\Phi)\\
(\text{pull back})=&\
D^{-1}\partial_z(p\circ\Phi)
-BD^{-1}\partial_{R_1}(p\circ\Phi)
-CD^{-1}\partial_{R_2}(p\circ\Phi).
\end{split}
\end{equation*}
Hence (omitting $\circ\Phi$) we obtain
\begin{equation*}
\begin{split}
-\partial_{R_1}(\nabla p\cdot \tau)\Big|_{R_1=R_2=0}
=&\
\Bigl((-\kappa+\partial_z\alpha_1-\alpha_1^2\kappa-T\alpha_2)\partial_z p
-\partial_{R_1}\partial_{z} p\\
&\quad+\alpha_1\kappa\,\partial_{R_1}p+T\,\partial_{R_2}p\Bigr)\Big|_{R_1=R_2=0}\\
(\text{commuting }\partial_{R_1}\text{ and }\partial_{z})=&\
(-\kappa+\partial_z\alpha_1-\alpha_1^2\kappa-T\alpha_2)\partial_z p
-\partial_z\partial_{R_1}p\\
&\quad+\alpha_1\kappa\,\partial_{R_1}p
+T\,\partial_{R_2}p\\
=:&\ (RHS).
\end{split}
\end{equation*}
\begin{remark}
We can rephrase the commutativity of $\partial_{R_1}$ and $\partial_z$ as
\[
[\partial_z,\partial_{R_1}]
:=\partial_z\partial_{R_1}-\partial_{R_1}\partial_z=0,
\]
where $[\cdot,\cdot]$ denotes the Lie bracket.
For previous studies exploiting this property, see Chan--Czubak--Y \cite[Section~2.5]{CCY},
Lichtenfelz--Y \cite{LY}, and Shimizu--Y \cite{SY}; for an earlier appearance, see Ma--Wang \cite[(3.7)]{MW}.
\end{remark}

Since (along $R_1=R_2=0$, again omitting $\circ\Phi$)
\begin{equation*}
\begin{split}
-\partial_{R_1}p
&=-(\nabla p\cdot n)-\alpha_1(\nabla p\cdot\tau)
=v^2\kappa+v'\alpha_1,\\
-\partial_z\partial_{R_1}p
&=\partial_z(v^2\kappa+v'\alpha_1)
=v^2\partial_z\kappa+2v'\kappa+v'\partial_z\alpha_1+\frac{v''}{v}\alpha_1,\\
-\partial_{R_2}p
&=-(\nabla p\cdot b)-\alpha_2(\nabla p\cdot \tau)
=v'\alpha_2,\\
-\partial_z\partial_{R_2}p
&=\partial_z(v'\alpha_2)
=\frac{v''}{v}\alpha_2+v'\partial_z\alpha_2,\\
-\partial_z p
&=-(\nabla p\cdot \tau)
=v',
\end{split}
\end{equation*}
we obtain
\begin{equation*}
\begin{split}
(RHS)
&=
-(-\kappa+\partial_z\alpha_1-\alpha_1^2\kappa-T\alpha_2)v'
+v^2\partial_z\kappa+2v'\kappa+v'\partial_z\alpha_1+\frac{v''}{v}\alpha_1\\
&\quad
-\alpha_1\kappa(v^2\kappa+v'\alpha_1)-Tv'\alpha_2\\
&=
v^2\partial_z\kappa+3v'\kappa+\frac{v''}{v}\alpha_1-v^2\kappa^2\alpha_1.
\end{split}
\end{equation*}

On the other hand, similarly,
\begin{equation*}
\begin{split}
-\partial_{R_2}(\nabla p\cdot \tau)
&=(\partial_z\alpha_2+T\alpha_1-\alpha_1\alpha_2\kappa)\partial_zp
-\partial_{R_2}\partial_{z} p-(T-\alpha_2\kappa)\partial_{R_1} p\\
(\text{commuting }\partial_{R_2}\text{ and }\partial_{z})&=
(\partial_z\alpha_2+T\alpha_1-\alpha_1\alpha_2\kappa)\partial_zp
-\partial_{z}\partial_{R_2} p-(T-\alpha_2\kappa)\partial_{R_1} p\\
&=
-(\partial_z\alpha_2+T\alpha_1-\alpha_1\alpha_2\kappa)v'
+\frac{v''}{v}\alpha_2+v'\partial_z\alpha_2\\
&\qquad+(T-\alpha_2\kappa)(v^2\kappa+v'\alpha_1)\\
&=
\frac{v''}{v}\alpha_2+v^2T\kappa-v^2\kappa^2\alpha_2,
\end{split}
\end{equation*}
for $R_1=R_2=0$.
Combining these identities with \eqref{pressure-estimate}, we arrive at the coupled closure equations
\begin{equation}\label{closure equation}
\begin{split}
-v'\kappa -v^2\partial_z\kappa+\partial_t^2\alpha_1-v^2\kappa^2\alpha_1
&=
v^2\partial_z\kappa+3v'\kappa+\frac{v''}{v}\alpha_1-v^2\kappa^2\alpha_1,\\
v^2T\kappa+\partial_t^2\alpha_2-v^2\kappa^2\alpha_2
&=
\frac{v''}{v}\alpha_2+v^2 T\kappa-v^2\kappa^2\alpha_2.
\end{split}
\end{equation}
Note that $\partial_z\kappa=v^{-1}\kappa'$ since $\partial_t=v\partial_z$ along the vortex axis.
Rearranging \eqref{closure equation} yields
\begin{equation*}
\begin{cases}
\displaystyle
\alpha_1''=\frac{v''}{v}\alpha_1+2v\kappa'+4v'\kappa,\\
\displaystyle
\alpha_2''=\frac{v''}{v}\alpha_2.
\end{cases}
\end{equation*}
These are the desired wave equations.

\section{Basic setup of numerical implementations: Bayesian optimization and adjoint method}

In this section, we prepare for the numerical implementation.
First we observe that the axis of swirling particles is expressed as 
\begin{equation}\label{zeta}
\zeta(t,s)=\tau(t,s)-\alpha_1(t,s)n(t,s)-\alpha_2(t,s)b(t,s).
\end{equation}
and we define the unit $s$-direction (vortex axis) by
\begin{equation}\label{es}
    \begin{split}
e_s(t,s):=
&\frac{\partial_s\Phi(t,s)}{|\partial_s\Phi(t,s)|}\\
=&
\bigl(e_s(t,s)\cdot\tau(t,s)\bigr)\tau(t,s)
+
\bigl(e_s(t,s)\cdot n(t,s)\bigr)n(t,s)
+
\bigl(e_s(t,s)\cdot b(t,s)\bigr)b(t,s).
\end{split}
\end{equation}
Throughout this paper we assume that the vortex axis and the axis of swirling particles
are perfectly aligned at the initial time.
That is, 
\begin{equation*}
    e_s(t,s)\cdot\frac{\zeta(t,s)}{|\zeta(t,s)|}\bigg|_{t=0}=1
    \quad\text{and}\quad
    \partial_t\left( e_s(t,s)\cdot\frac{\zeta(t,s)}{|\zeta(t,s)|}\right)\bigg|_{t=0}=0,
\end{equation*}
in other words,
\begin{equation*}
    (e_s(t,s)\cdot \tau(t,s))\zeta(t,s)=e_s(t,s),\quad
    \partial_t\left((e_s(t,s)\cdot \tau(t,s))\zeta(t,s)\right)=\partial_te_s(t,s)
\end{equation*}
at $t=0$.
Comparing the coefficients in the Frenet frame yields
\begin{equation*}
\begin{split}
\alpha_1(0,s)
=
-\frac{e_s(0,s)\cdot n(0,s)}{e_s(0,s)\cdot\tau(0,s)},
\quad
&\alpha_2(0,s)
=
-\frac{e_s(0,s)\cdot b(0,s)}{e_s(0,s)\cdot\tau(0,s)},\\
\partial_t\alpha_1(0,s)
=
-\partial_t\left(\frac{e_s(t,s)\cdot n(t,s)}{e_s(t,s)\cdot\tau(t,s)}\right)\bigg|_{t=0},
\quad
&\partial_t\alpha_2(0,s)
=
-\partial_t\left(\frac{e_s(t,s)\cdot b(t,s)}{e_s(t,s)\cdot\tau(t,s)}\right)\bigg|_{t=0}.
\end{split}
\end{equation*}
If the two axes are not parallel, the resulting rotating flow tends to exhibit shear.
From a physical viewpoint, it is therefore natural to assume that the system favors configurations in which these axes remain aligned during the evolution, thereby minimizing shear.
Motivated by this consideration, we introduce the following
mean square directional correlation
(MSDC) as a functional:
\begin{equation}\label{MSDC}
J(\Phi):=
\frac{1}{t_1-t_0}
\int_{t_0}^{t_1}
\int_0^1
\left|
e_s(t,s)
\cdot
\frac{\zeta(t,s)}{|\zeta(t,s)|}
\right|^2
\, ds\, dt,
\end{equation}
which we adopt as an objective functional to be maximized.

\subsection{Bayesian optimization}

We first prescribe a transport flow that expands in the radial direction,
characterized by a scalar function $\Gamma(t)$ satisfying
\[
\Gamma(t_0)=0, \qquad \Gamma(t)>0, \qquad \Gamma'(t)>0
\quad\text{for}\quad t\in (t_0,t_1).
\]
The time evolution of the ring shape, for $t \in [t_0,t_1)$ and $s \in [0,1)$,
is defined by
\begin{equation*}
\Phi(t,s)
=
\epsilon K(t,s)+\Bigl(1+\Gamma(t)+\gamma_1(t,s)
\Bigr)
\begin{pmatrix}
\cos 2\pi s\\
\sin 2\pi s\\
0
\end{pmatrix}
+
\gamma_2(t,s)
\begin{pmatrix}
0\\
0\\
1
\end{pmatrix}.
\end{equation*}
Here, $K$ is a Kelvin wave defined by
\begin{equation*}
K(t,s)
=
\begin{pmatrix}
\cos\big(2\pi m s - \omega t\big)\cos(2\pi s) \\
\cos\big(2\pi m s - \omega t\big)\sin(2\pi s) \\
\sin\big(2\pi m s - \omega t\big)
\end{pmatrix}.
\end{equation*}
The functions $\gamma_1(t,s)$ and $\gamma_2(t,s)$ describe a time-dependent wavy deformation and are given by
\begin{equation*}
\begin{split}
\gamma_{1}(t,s)
&=\sum_{j=0}^{J}\sum_{k=0}^{K}
\Bigl(
c^{11}_{jk}\left(\frac{t-t_0}{t_1-t_0}\right)^{j+1}\sin 2\pi ks
+
c^{12}_{jk}\left(\frac{t-t_0}{t_1-t_0}\right)^{j+1}\cos 2\pi ks
\Bigr),\\
\gamma_{2}(t,s)
&=\sum_{j=0}^{J}\sum_{k=0}^{K}
\Bigl(
c^{21}_{jk}\left(\frac{t-t_0}{t_1-t_0}\right)^{j+1}\sin 2\pi ks
+
c^{22}_{jk}\left(\frac{t-t_0}{t_1-t_0}\right)^{j+1}\cos 2\pi ks
\Bigr),
\end{split}
\end{equation*}
where the coefficients $\{c^{\ell m}_{jk}\}_{jk\ell m}\subset\mathbb{R}$ are treated as learnable parameters for the Bayesian optimization.
The hyperparameters are $\epsilon$, $m$, $\omega$, $J$ and  $K$.

\subsection{Adjoint-based optimization}
After obtaining an optimized vortex-ring evolution via Bayesian optimization as a global search step, we then apply adjoint-based optimization.
In what follows, we present a step-by-step formulation of the adjoint-based optimization procedure.
Recall that 
\[
\partial_t \Phi = v \tau,
\quad
v = |\partial_t \Phi|,
\quad
\tau = \frac{\partial_t \Phi}{|\partial_t \Phi|}
\quad\text{and}\quad
\partial_t^2 \Phi = v' \tau + v^2 \kappa n.
\]
Therefore, the curvature can be written as
\begin{equation*}
\begin{split}
\kappa
=&
\frac{|A|}{v^2},\qquad A:=\partial_t^2\Phi-\frac{v'}{v}\partial_t\Phi= v^2\kappa n.
\end{split}
\end{equation*}

\begin{remark}
In the continuous setting, the identity
\begin{equation*}
\partial_t | \partial_t \Phi |
= \frac{\partial_t \Phi}{|\partial_t \Phi|} \cdot \partial_{tt} \Phi
\end{equation*}
follows from the chain rule, which implies that the vector
\begin{equation*}
A = \partial_{tt} \Phi - \frac{\partial_t |\partial_t \Phi|}{|\partial_t \Phi|} \,\partial_t \Phi
\end{equation*}
is orthogonal to the tangent direction, i.e.,\ $A \cdot \tau = 0$.
However, after time discretization, the discrete derivative $D_t$ does not satisfy the chain rule. In particular
\begin{equation*}
D_t |\partial_t \Phi| \neq \frac{\partial_t \Phi}{|\partial_t \Phi|} \cdot D_t (\partial_t \Phi).
\end{equation*}
As a consequence, the discrete analogue of $A$ is no longer exactly orthogonal to $\tau$, and a spurious tangential component may appear. For numerical stability, we therefore enforce orthogonality by projection, replacing $A$ with
\begin{equation*}
A_{\perp} = A - (A \cdot \tau)\,\tau.
\end{equation*}
\end{remark}
Direct computations yield
\begin{equation}
\label{delta v}
\begin{split}\delta v
&=
\frac{\partial_t\Phi}{|\partial_t\Phi|}\cdot \partial_t\delta\Phi
=
\tau\cdot\partial_t\delta\Phi\\
\delta\tau&=
\delta\left(\frac{1}{v}\right)\partial_t\Phi+\frac{\partial_t\delta\Phi}{v}
=\frac{1}{v}(I-\tau\otimes\tau)\partial_t\delta\Phi\\
\end{split}
\end{equation}
and then
\begin{equation}\label{delta kappa}
\delta\kappa
=
-\frac{2\kappa}{v}\delta v
+\frac{1}{v^2}\delta{|A|}
\quad\text{and}\quad 
\delta{|A|}
=
\frac{A}{|A|}\cdot \delta A
=
n\cdot \delta A.
\end{equation}
On the other hand,
\[
\delta A
=
\partial_t^2\delta\Phi
-
\delta{\left(\frac{v'}{v}\right)}\partial_t\Phi
-
\frac{v'}{v}\partial_t\delta \Phi
\quad\text{and}\quad
\delta {\left(\frac{v'}{v}\right)}
=
\partial_t\left(\frac{\delta v}{v}\right).
\]
Therefore
\begin{equation}\label{delta A}
\delta A
=
\partial_t^2\delta\Phi
-
\partial_t\left(\frac{\delta v}{v}\right)\partial_t\Phi
-
\frac{v'}{v}\partial_t\delta\Phi.
\end{equation}
Now we formulate the variational formulas of $\kappa$ and $n$.
Substituting this into \eqref{delta kappa} and since \(n\cdot\partial_t\Phi=0\), we obtain
\begin{equation}
\delta\kappa
=
\frac1{v^2}n\cdot \partial_t^2\delta\Phi
-\frac{v'}{v^3}n\cdot \partial_t\delta\Phi
-\frac{2\kappa}{v}\,\tau\cdot\partial_t\delta\Phi.
\end{equation}
Recall that $\displaystyle n=\frac{A}{|A|}$.
By \eqref{delta A}, we obtain
\begin{equation*}
    \begin{split}
\delta n
=&
\frac1{|A|}
\left(I-n\otimes n\right)\delta A\\
=&
\frac1{v^2\kappa}
\left(I-n\otimes n\right)
\left(
\partial_t^2\delta\Phi
-
\partial_t\left(\frac{\delta v}{v}\right)\partial_t\Phi
-
\frac{v'}{v}\partial_t\delta\Phi
\right).
\end{split}
\end{equation*}
Since \(n\cdot\partial_t\Phi=0\),
we have 
\begin{equation*}
 \begin{split}
\delta n
=
&
\frac1{v^2\kappa}
\left(I-n\otimes n\right)\partial_t^2\delta\Phi
-
\frac1{v^2\kappa}
\partial_t\left(\frac{\delta v}{v}\right)\partial_t\Phi
-
\frac{v'}{v^3\kappa}
\left(I-n\otimes n\right)\partial_t\delta\Phi\\
=:& I_1+I_2+I_3.
\end{split}
\end{equation*}
Using \eqref{delta v},
\begin{equation*}
(I-n\otimes n)X = (\tau\cdot X)\tau + (b\cdot X)b
\quad
\text{and}\quad
\partial_t\tau = v\kappa n,
\end{equation*}
we have 
\begin{equation*}
\partial_t\!\left(\frac{\delta v}{v}\right)
=
\kappa n\cdot \partial_t\delta\Phi
+\frac{1}{v}\tau\cdot \partial_t^2\delta\Phi
-\frac{v'}{v^2}\tau\cdot \partial_t\delta\Phi.
\end{equation*}
Thus
\begin{equation*}
I_2=
-\frac{n\cdot \partial_t\delta\Phi}{v}\tau
-\frac{\tau\cdot \partial_t^2\delta\Phi}{v^2\kappa}\tau
+\frac{v'}{v^3\kappa}(\tau\cdot \partial_t\delta\Phi)\tau.
\end{equation*}
On the other hand, we see 
\begin{equation*}
I_1=
\frac{\tau\cdot \partial_t^2\delta\Phi}{v^2\kappa}\tau
+
\frac{b\cdot \partial_t^2\delta\Phi}{v^2\kappa}b
\quad\text{and}\quad 
I_3=
-\frac{v'}{v^3\kappa}(\tau\cdot \partial_t\delta\Phi)\tau
-\frac{v'}{v^3\kappa}(b\cdot \partial_t\delta\Phi)b.
\end{equation*}
This means that 
the $\tau$-components cancel, and we obtain
\begin{equation*}
\delta n
=
-\frac{n\cdot \partial_t\delta\Phi}{v}\,\tau
+
\left(
\frac{b\cdot \partial_t^2\delta\Phi}{v^2\kappa}
-\frac{v'}{v^3\kappa}\, b\cdot \partial_t\delta\Phi
\right)b.
\end{equation*}
Now let us derive $\delta b$. Using the Frenet relation \(b=\tau\times n\), we see that
\[
b
=
\frac{1}{v^3\kappa}
\,\partial_t \Phi \times \partial_t^2 \Phi .
\]
We formulate the variational formula of this $b$. By 
\[
\begin{cases}
\delta(\partial_t\Phi\times\partial_t^2\Phi)
=
\partial_t\delta\Phi\times\partial_t^2\Phi
+
\partial_t\Phi\times\partial_t^2\delta\Phi,\\
\displaystyle\delta{\left(\frac{1}{v^3\kappa}\right)}
=
-\frac{1}{v^3\kappa}
\left(
3\frac{\delta v}{v}
+
\frac{\delta\kappa}{\kappa}
\right),
\end{cases}
\]
we have 
\[
\delta b
=
\frac{1}{v^3\kappa}
\left(
\partial_t\delta\Phi\times\partial_t^2\Phi
+
\partial_t\Phi\times\partial_t^2\delta\Phi
\right)
-
\left(
3\frac{\delta v}{v}
+
\frac{\delta\kappa}{\kappa}
\right)b.
\]
Therefore
\begin{align*}
\delta b
&=
\frac{1}{v^3\kappa}
\Bigl[
v'\bigl(\partial_t\delta\Phi\times\tau\bigr)
+
v^2\kappa\bigl(\partial_t\delta\Phi\times n\bigr)
+
v\bigl(\tau\times\partial_t^2\delta\Phi\bigr)
\Bigr]
-
\left(
3\frac{\delta v}{v}
+
\frac{\delta\kappa}{\kappa}
\right)b \\
&=
\frac{1}{v^3\kappa}
\Bigl[
v'\bigl(
-(n\cdot\partial_t\delta\Phi)b
+
(b\cdot\partial_t\delta\Phi)n
\bigr)
+
v^2\kappa\bigl(
(\tau\cdot\partial_t\delta\Phi)b
-
(b\cdot\partial_t\delta\Phi)\tau
\bigr) \\
&\qquad\qquad
+
v\bigl(
(n\cdot\partial_t^2\delta\Phi)b
-
(b\cdot\partial_t^2\delta\Phi)n
\bigr)
\Bigr]
-
\left(
3\frac{\delta v}{v}
+
\frac{\delta\kappa}{\kappa}
\right)b \\
&=
-\frac{b\cdot\partial_t\delta\Phi}{v}\,\tau
+
\left(
\frac{v'}{v^3\kappa}\,b\cdot\partial_t\delta\Phi
-
\frac{1}{v^2\kappa}\,b\cdot\partial_t^2\delta\Phi
\right)n \\
&\quad
+
\Bigl(
-\frac{v'}{v^3\kappa}\,n\cdot\partial_t\delta\Phi
+
\frac{1}{v^2\kappa}\,n\cdot\partial_t^2\delta\Phi
+
\frac{1}{v}\,\tau\cdot\partial_t\delta\Phi
\Bigr)b
-
\left(
3\frac{\delta v}{v}
+
\frac{\delta\kappa}{\kappa}
\right)b \\
&=
-\frac{b\cdot\partial_t\delta\Phi}{v}\,\tau
-
\left(
\frac{b\cdot\partial_t^2\delta\Phi}{v^2\kappa}
-
\frac{v'}{v^3\kappa}\,b\cdot\partial_t\delta\Phi
\right)n.
\end{align*}
Recall that
\[
e_s(t,s):=\frac{\partial_s\Phi(t,s)}{|\partial_s\Phi(t,s)|}.
\]
Finally we formulate the variational formula of $e_s$, by using the standard variation formula for a normalized vector: 
\begin{equation*}
\begin{split}
\delta e_s
=
\frac{1}{|\partial_s\Phi|}
\left(I-e_s\otimes e_s\right)\partial_s\delta\Phi
=
\frac{1}{|\partial_s\Phi|}
\left[
\partial_s\delta\Phi
-
\bigl(e_s\cdot \partial_s\delta\Phi\bigr)e_s
\right].
\end{split}
\end{equation*}

\subsection{Variational formulas}

Here we summarize the variational formulas:
\begin{equation}\label{formulas}
\begin{split}
\delta v
&=\tau\cdot\partial_t\delta\Phi,\\
\delta\tau
&=
\frac1v(I-\tau\otimes\tau)\partial_t\delta\Phi
=
\frac{n\cdot \partial_t\delta\Phi}{v}\,n
+
\frac{b\cdot \partial_t\delta\Phi}{v}\,b,\\
\delta n
&=
-\frac{n\cdot\partial_t\delta\Phi}{v}\,\tau
+
\left(
\frac{b\cdot\partial_t^2\delta\Phi}{v^2\kappa}
-\frac{v'}{v^3\kappa}\,b\cdot\partial_t\delta\Phi
\right)b,\\
\delta b
&=
-\frac{b\cdot\partial_t\delta\Phi}{v}\,\tau
-
\left(
\frac{b\cdot\partial_t^2\delta\Phi}{v^2\kappa}
-\frac{v'}{v^3\kappa}\,b\cdot\partial_t\delta\Phi
\right)n,\\
\delta\kappa
&=
\frac1{v^2}n\cdot\partial_t^2\delta\Phi
-\frac{v'}{v^3}n\cdot\partial_t\delta\Phi
-\frac{2\kappa}{v}\tau\cdot\partial_t\delta\Phi,\\
\delta e_s
&=
\frac{1}{|\partial_s\Phi|}
\left(I-e_s\otimes e_s\right)\partial_s\delta\Phi.\\
\end{split}
\end{equation}

\subsection{Evaluation function}

Recall the following mean square directional correlation (MSDC):
\begin{equation*}
J(\Phi)
:=
\frac{1}{t_1-t_0}
\int_{t_0}^{t_1}
\int_0^1
\left|
e_s(t,s)
\cdot
\frac{\zeta(t,s)}{|\zeta(t,s)|}
\right|^2
\, ds\, dt,
\end{equation*}
which serves as the objective functional to be maximized.
For simplicity, we write
\[
J(\Phi)
=
\frac{1}{t_1-t_0}
\int_{t_0}^{t_1}\int_0^1 c^2\,ds\,dt
\quad\text{for}\quad c:=e_s\cdot \hat\zeta
\quad\text{and}\quad\hat\zeta:=\frac{\zeta}{|\zeta|}
.
\]
Consider the following variation:
\[
\delta J(\delta\Phi)
=
\frac{2}{t_1-t_0}
\int_{t_0}^{t_1}\int_0^1
c\,\delta c
\,ds\,dt,\quad \delta c
=
\delta e_s\cdot \hat\zeta
+
e_s\cdot \delta {\hat\zeta}
\]
subject to the following initial and terminal conditions:
\begin{equation*}
    \delta\Phi(t_0)=\delta\Phi(t_1)=0\quad\text{and}\quad 
    \partial_t\delta\Phi(t_0)=\partial_t\delta\Phi(t_1)=0.
\end{equation*}
Since
\[
\delta{\hat\zeta}
=
\frac1{|\zeta|}
\left(I-\hat\zeta\otimes\hat\zeta\right)\delta \zeta,
\]
we have
\[
e_s\cdot \delta {\hat\zeta}
=
\frac1{|\zeta|}
\bigl(e_s-(e_s\cdot\hat\zeta)\hat\zeta\bigr)\cdot\delta\zeta
=
\frac1{|\zeta|}(e_s-c\hat\zeta)\cdot\delta\zeta.
\]
Moreover, 
\[
\delta e_s\cdot\hat\zeta
=
\frac1{|\partial_s\Phi|}
(\hat\zeta-c\,e_s)\cdot\partial_s\delta\Phi.
\]
Hence
\begin{equation}\label{delta MSDC}
\begin{split}
\delta J(\delta \Phi)
&=
\frac{2}{t_1-t_0}
\int_{t_0}^{t_1}\int_0^1
c\left[
\frac1{|\partial_s\Phi|}
(\hat\zeta-c\,e_s)\cdot\partial_s\delta\Phi
+
\frac1{|\zeta|}
(e_s-c\hat\zeta)\cdot\delta\zeta
\right]
\,ds\,dt\\
&=
\frac{2}{t_1-t_0}
\int_{t_0}^{t_1}\int_0^1
c\bigg[
\frac1{|\partial_s\Phi|}
(\hat\zeta-c\,e_s)\cdot\partial_s\delta\Phi\\
&\qquad +
\frac1{|\zeta|}
(e_s-c\hat\zeta)\cdot
\bigl(
\delta\tau-\delta\alpha_1 n-\alpha_1\delta n-\delta\alpha_2 b-\alpha_2\delta b
\bigr)
\bigg]
\,ds\,dt\\
&
=
\int_{t_0}^{t_1}\int_0^1
\left[a\cdot \partial_s\delta\Phi
+
m\cdot(\delta\tau-\alpha_1\delta n-\alpha_2\delta b)
+
q_1\,\delta\alpha_1
+
q_2\,\delta\alpha_2\right]
\,ds\,dt,
\end{split}
\end{equation}
with
\[
a
:=
\frac{2}{t_1-t_0}\,
\frac{c}{|\partial_s\Phi|}
(\hat\zeta-c\,e_s),
\qquad
m
:=
\frac{2}{t_1-t_0}\,
\frac{c}{|\zeta|}
(e_s-c\,\hat\zeta),
\]
\[
q_1:=-m\cdot n,
\qquad
q_2:=-m\cdot b.
\]
The main objective is to find an explicit representation of $\tilde S$ such that  
\begin{equation*}
\iint q_1\,\delta\alpha_1
+
q_2\,\delta\alpha_2
=\iint \tilde S\cdot \delta\Phi.
\end{equation*}

\subsection{Variational formulas for $\alpha_1$ and $\alpha_2$ and adjoint representations.}
Recall that $\alpha_1$ and $\alpha_2$
satisfy the following wave-type equations:
\begin{equation*}
\alpha_1''
=
\frac{v''}{v}\alpha_1
+2v\kappa'
+4v'\kappa,
\qquad
\alpha_2''
=
\frac{v''}{v}\alpha_2.
\end{equation*}
Then we obtain
\begin{equation*}
    \begin{cases}\displaystyle
\delta\alpha_1''
=
\delta{\left(\frac{v''}{v}\right)}\alpha_1
+
\frac{v''}{v}\delta\alpha_1
+
2\delta v\,\kappa'
+
2v\,\delta\kappa'
+
4\delta v'\kappa
+
4v'\delta\kappa,\\
\displaystyle
\delta\alpha_2''
=
\delta{\left(\frac{v''}{v}\right)}\alpha_2
+
\frac{v''}{v}\delta\alpha_2.
\end{cases}
\end{equation*}
Since
\[
\delta{\left(\frac{v''}{v}\right)}
=
\frac{\delta v''}{v}-\frac{v''}{v^2}\delta v,
\]
it follows that
\[
\delta\alpha_1''-\frac{v''}{v}\delta\alpha_1=S_1[\delta\Phi],
\qquad
\delta\alpha_2''-\frac{v''}{v}\delta\alpha_2=S_2[\delta\Phi],
\]
with
\begin{equation}\label{S}
\begin{cases}
\displaystyle
S_1[\delta\Phi]
=
\left(\frac{\delta v''}{v}-\frac{v''}{v^2}\delta v\right)\alpha_1
+2\,\delta v\,\kappa'
+2v\,\delta\kappa'
+4\,\delta v'\,\kappa
+4v'\,\delta\kappa,\\
\displaystyle
S_2[\delta\Phi]
=
\left(\frac{\delta v''}{v}-\frac{v''}{v^2}\delta v\right)\alpha_2.
\end{cases}
\end{equation}
We aim
to rewrite $\delta\alpha_1$ and $\delta\alpha_2$ in \eqref{delta MSDC} in terms of  $\delta\Phi$ explicitly.
To this end, we introduce adjoint variables \(\lambda_1,\lambda_2\) solving
\[
\lambda_1''-\frac{v''}{v}\lambda_1=q_1,
\qquad
\lambda_2''-\frac{v''}{v}\lambda_2=q_2
\]
with the terminal conditions
\[
\lambda_j(t_1,s)=0,
\qquad
\partial_t\lambda_j(t_1,s)=0. 
\]
On the other hand the initial conditions are given by 
\[
\delta\alpha_j(t_0,s)=0,
\qquad
\partial_t\delta\alpha_j(t_0,s)=0.
\]
\begin{remark}
These initial conditions are justified by assuming 
\begin{equation*}
\delta\Phi(t_0)=\partial_t\delta\Phi(t_0)=\partial_t^2\delta\Phi(t_0)=\partial_t^3\delta\Phi(t_0)=0.
\end{equation*}
This will be explained later.
\end{remark}
Then we obtain 
\begin{equation}\label{qalpha}
\iint q_1\,\delta\alpha_1+q_2\,\delta\alpha_2
=
\iint \lambda_1 S_1[\delta\Phi]+\lambda_2 S_2[\delta\Phi].
\end{equation}
After integration by parts in $t$, we obtain
\[
\iint
\lambda_1 S_1[\delta\Phi]+\lambda_2 S_2[\delta\Phi]
=
\iint
\Lambda_v\,\delta v
+
\Lambda_\kappa\,\delta\kappa,
\]
where
\[
\Lambda_v
=
\left(\frac{\eta}{v}\right)''
-4(\kappa\lambda_1)'
+2\kappa'\lambda_1
-\frac{v''}{v^2}\eta,\quad
\eta:=\alpha_1\lambda_1+\alpha_2\lambda_2
\]
and
\[
\Lambda_\kappa
=
-(2v\lambda_1)'
+4v'\lambda_1
=
2v'\lambda_1-2v\lambda_1'.
\]
Write
\[
m_\tau:=m\cdot\tau,
\qquad
m_n:=m\cdot n,
\qquad
m_b:=m\cdot b.
\]
Then from \eqref{delta MSDC} and the formulas \eqref{formulas}, we finally obtain
\[
\delta J(\delta \Phi)
=
\int_{t_0}^{t_1}\int_0^1
\left(a\cdot \partial_s\delta\Phi
+
B\cdot \partial_t\delta\Phi
+
C\cdot \partial_t^2\delta\Phi\right)
\,ds\,dt,
\]
where
\begin{equation*}
\begin{split}
B
&=
\left(\Lambda_v-\frac{2\kappa}{v}\Lambda_\kappa\right)\tau
+
\left(
\frac{m_n+\alpha_1 m_\tau}{v}
-\frac{v'}{v^3}\Lambda_\kappa
\right)n\\
&\qquad
+
\left(
\frac{m_b+\alpha_2 m_\tau}{v}
+\frac{v'}{v^3\kappa}(\alpha_1 m_b-\alpha_2 m_n)
\right)b,
\end{split}
\end{equation*}
and
\[
C
=
\frac{\Lambda_\kappa}{v^2}\,n
+
\frac{-\alpha_1 m_b+\alpha_2 m_n}{v^2\kappa}\,b.
\]
Integrating by parts, we finally obtain
\[
\delta J(\delta \Phi)
=
\int_{t_0}^{t_1}\int_0^1
\left(
-\partial_s a
-\partial_t B
+\partial_t^2 C
\right)\cdot \delta\Phi
\,ds\,dt.
\]

\subsection{Variational formulas for the initial conditions of $\delta\alpha_1$ and $\delta\alpha_2$}
Recall the initial condition:
\begin{equation*}
\begin{split}
\alpha_1(0,s)
=
-\frac{e_s(0,s)\cdot n(0,s)}{e_s(0,s)\cdot\tau(0,s)},
\quad
&\alpha_2(0,s)
=
-\frac{e_s(0,s)\cdot b(0,s)}{e_s(0,s)\cdot\tau(0,s)},\\
\partial_t\alpha_1(0,s)
=
-\partial_t\left(\frac{e_s(t,s)\cdot n(t,s)}{e_s(t,s)\cdot\tau(t,s)}\right)\bigg|_{t=0},
\quad
&\partial_t\alpha_2(0,s)
=
-\partial_t\left(\frac{e_s(t,s)\cdot b(t,s)}{e_s(t,s)\cdot\tau(t,s)}\right)\bigg|_{t=0}.
\end{split}
\end{equation*}
For \(\alpha_1\), we have
\[
\delta\alpha_1
=
-\frac{
\bigl(\delta e_s\cdot n+e_s\cdot\delta n\bigr)(e_s\cdot\tau)
-
(e_s\cdot n)\bigl(\delta e_s\cdot\tau+e_s\cdot\delta\tau\bigr)
}{
(e_s\cdot\tau)^2
}.
\]
Similarly, for \(\alpha_2\), we have 
\[
\delta\alpha_2
=
-\frac{
\bigl(\delta e_s\cdot b+e_s\cdot\delta b\bigr)(e_s\cdot\tau)
-
(e_s\cdot b)\bigl(\delta e_s\cdot\tau+e_s\cdot\delta\tau\bigr)
}{
(e_s\cdot\tau)^2
}.
\]
Thus, under the assumption that
\[
\delta\Phi(t_0)=\partial_t\delta\Phi(t_0)=\partial_t^2\delta\Phi(t_0)=\partial_t^3\delta\Phi(t_0)=0,
\]
and using \eqref{formulas}, we obtain
\[
\delta\alpha_1(t_0)=\delta\alpha_2(t_0)=0,
\quad
\partial_t\delta\alpha_1(t_0)=\partial_t\delta\alpha_2(t_0)=0.
\]

\subsection{Discrete formulation}
In the previous subsection, we presented the continuous adjoint formulation. 
In this subsection, we introduce its discrete counterpart. 
Although the essential idea remains the same, we summarize the discrete procedure 
to facilitate understanding of the implementation. 
The key ingredient is the use of backpropagation.
For
\begin{equation*}
\Phi=\{\Phi_j\}_{j=0}^{N},
\qquad
\Phi_j\in\mathbb{R}^3,
\end{equation*}
let
\begin{equation*}
\alpha=\{\alpha_j\}_{j=0}^{N},
\qquad
\alpha_j=
\begin{pmatrix}
\alpha_{1,j}\\
\alpha_{2,j}
\end{pmatrix}
\in\mathbb{R}^2,
\end{equation*}
be defined by the discrete evolution (corresponding to \eqref{wave equations})
\begin{equation}\label{eq:discrete_alpha_update}
\alpha_{j+1}
=
F\bigl(\alpha_j,\alpha_{j-1},\{\Phi_{j+\ell}\}_{\ell\in\Lambda}\bigr),
\qquad
j=1,\dots,N-1,
\end{equation}
where $\Lambda\subset\mathbb Z$ is a finite index set. 
This finite-stencil representation encodes the dependence on $\Phi$ through neighboring time levels and allows us to treat, in a unified way, the higher-order time derivatives appearing in the continuous formulation. 
The contribution of $\alpha$ to the objective functional is written as
\begin{equation*}
\delta J_\alpha(\delta\Phi)
:=
\sum_{j=0}^{N-2}
\langle q_j,\delta\alpha_j\rangle,
\qquad
q_j=
\begin{pmatrix}
q_{1,j}\\
q_{2,j}
\end{pmatrix},
\end{equation*}
which corresponds to \eqref{qalpha}.
Taking the variation of \eqref{eq:discrete_alpha_update}, we obtain the linearized forward equation 
\begin{equation}\label{eq:discrete_linearized_alpha}
\delta\alpha_{j+1}
=
A_j\,\delta\alpha_j
+
B_j\,\delta\alpha_{j-1}
+
r_j,
\qquad
j=1,\dots,N-1,
\end{equation}
where
\begin{equation*}
A_j:=\partial_{\alpha_j}F,
\qquad
B_j:=\partial_{\alpha_{j-1}}F,
\qquad
r_j
=
\sum_{\ell\in\Lambda}
S_{j,\ell}\,\delta\Phi_{j+\ell},
\end{equation*}
and each $S_{j,\ell}$ is a linear operator acting on $\delta\Phi_{j+\ell}$ (corresponding to \eqref{S}).
We impose endpoint conditions on the  variations of $\alpha$, namely,
\begin{equation*}
\delta\alpha_0=0,
\qquad
\delta\alpha_1=0.
\end{equation*}
Now we introduce the adjoint variable $\lambda_j\in\mathbb{R}^2$ and add the identity
\begin{equation*}
0
=
\sum_{j=1}^{N-1}
\left\langle
\lambda_{j+1},
\delta\alpha_{j+1}
-
A_j\delta\alpha_j
-
B_j\delta\alpha_{j-1}
-
r_j
\right\rangle,
\end{equation*}
which holds by \eqref{eq:discrete_linearized_alpha}. 
Hence,
\begin{equation*}
\delta J_\alpha
=
\sum_{j=0}^{N-2}
\langle q_j,\delta\alpha_j\rangle
+
\sum_{j=1}^{N-1}
\left\langle
\lambda_{j+1},
\delta\alpha_{j+1}
-
A_j\delta\alpha_j
-
B_j\delta\alpha_{j-1}
-
r_j
\right\rangle.
\end{equation*}
Rearranging the terms with respect to $\delta\alpha_j$, we obtain
\begin{equation*}
\delta J_\alpha
=
\sum_{j=2}^{N-2}
\left\langle
q_j+\lambda_j-A_j^\top\lambda_{j+1}-B_{j+1}^\top\lambda_{j+2},
\delta\alpha_j
\right\rangle
-
\sum_{j=1}^{N-1}
\langle \lambda_{j+1},r_j\rangle.
\end{equation*}
We therefore define the backward recursion
\begin{equation}\label{eq:discrete_adjoint_alpha}
\lambda_j
=
A_j^\top\lambda_{j+1}
+
B_{j+1}^\top\lambda_{j+2}
-
q_j,
\qquad
j=N-2,\dots,1,
\end{equation}
together with the terminal conditions
\begin{equation*}
\lambda_{N}=0,
\qquad
\lambda_{N-1}=0.
\end{equation*}
Then all terms containing $\delta\alpha_j$ vanish, and we are left with
\begin{equation*}
\delta J_\alpha
=
-
\sum_{j=1}^{N-1}
\langle \lambda_{j+1},r_j\rangle.
\end{equation*}
Substituting the finite-stencil form of $r_j$, we obtain
\begin{equation*}
\delta J_\alpha
=
-
\sum_{j=1}^{N-1}
\sum_{\ell\in\Lambda}
\left\langle
\lambda_{j+1},
S_{j,\ell}\,\delta\Phi_{j+\ell}
\right\rangle.
\end{equation*}
Taking transposes and grouping the coefficients of each $\delta\Phi_m$, we finally obtain
\begin{equation*}
\delta J_\alpha
=
\sum_{m=0}^{N}
\left\langle
\tilde S_m,\delta\Phi_m
\right\rangle,\quad
\tilde S_m
:=
-
\sum_{\substack{j=1,\dots,N-1\\ j+\ell=m}}
S_{j,\ell}^\top\lambda_{j+1},
\end{equation*}
where we implicitly assume boundary conditions for $\delta\Phi$, so that contributions outside the index range vanish.

\begin{remark}
The precise forward evolution of $\alpha_j$ is given by 
\begin{equation*}
\alpha_{j+1}
=
(2 + \Delta t^2 a_j)\,\alpha_j
-
\alpha_{j-1}
+
\Delta t^2 f_j,\quad a_j=\frac{v''_j}{v_j}\quad\text{and}\quad
f_j=
\begin{pmatrix}
2v_j\kappa_j'+4v_j'\kappa_j\\
0
\end{pmatrix}
\end{equation*}
which corresponds to \eqref{eq:discrete_alpha_update}.
Therefore, the linearized forward equation takes the form
\begin{equation*}
\delta \alpha_{j+1}
=
A_j \, \delta \alpha_j
+
B_j \, \delta \alpha_{j-1}
+
r_j,
\end{equation*}
with
\begin{equation*}
A_j = 2 + \Delta t^2 a_j,
\quad
B_j = -1\quad
\text{and}
\quad
r_j
=
\Delta t^2 \alpha_j \, \delta a_j
+
\Delta t^2 \, \delta f_j.
\end{equation*}
Then the discrete adjoint equation  is given by \eqref{eq:discrete_adjoint_alpha}.
\end{remark}

\subsection{Projected gradient under temporal admissibility and spatial band limitation}

Let $G:=\nabla_{\Phi}J$ denote the formal $L^2$-gradient obtained from the previous subsections.
The actual update direction must satisfy the temporal admissibility conditions at the initial and terminal times, and in addition we may restrict the perturbation to a finite spatial Fourier band in the parameter $s$.
First, we construct a temporal projection. We seek an admissible perturbation $\delta\Phi$ satisfying
\begin{equation}\label{admissible_condition}
\begin{split}
&\delta\Phi(t_0,s)
=\delta\Phi(t_1,s)=0,
\qquad
\partial_t\delta\Phi(t_0,s)=\partial_t\delta\Phi(t_1,s)=0,\\
&    \partial_t^2\delta\Phi(t_0,s)=0,\ 
    \partial_t^3\delta\Phi(t_0,s)=0
    \quad (\text{required for}\  \delta\alpha_j(t_0)=\partial_t\delta\alpha_j(t_0)=0)
    \end{split}
\end{equation}
which maximizes the first variation
\[
\delta J(\delta\Phi)=\int_{t_0}^{t_1}\int_0^1 G\cdot \delta\Phi \,ds\,dt.
\]
To construct 
the projection
in a spectral manner, we introduce the normalized time variable
\[
\xi:=\frac{t-t_0}{t_1-t_0},\qquad 
\xi\in[0,1],
\]
and write $\widetilde G(\xi,s):=G(t_0+(t_1-t_0)\xi,s)$.
We consider the boundary-adapted Fourier basis.
Let
\[
w(\xi):=\xi^4(1-\xi)^2,\qquad \xi\in[0,1],
\]
and define
\begin{align*}
\phi_0(\xi) &:= w(\xi),\\
\phi_k^{(c)}(\xi) &:= w(\xi)\cos(2\pi k\xi),\qquad k\ge1,\\
\phi_k^{(s)}(\xi) &:= w(\xi)\sin(2\pi k\xi),\qquad k\ge1.
\end{align*}
These functions satisfy
\[
\phi(0)=\phi(1)=0,\ \phi'(0)=\phi'(1)=0,\ \phi''(0)=0,\ \phi'''(0)=0
\]
and hence automatically fulfill the admissibility conditions.
Let
\[
V_N:=\mathrm{span}\left\{
\phi_0,\ \phi_k^{(c)},\ \phi_k^{(s)}\ :\ 1\le k\le N
\right\}.
\]
We define the projection
\[
P_t^{(N)} G(\xi,s)
=
a_0(s)\phi_0(\xi)
+
\sum_{k=1}^N
\Bigl(
a_k(s)\phi_k^{(c)}(\xi)
+
b_k(s)\phi_k^{(s)}(\xi)
\Bigr),
\]
as the $L^2$-orthogonal projection of $\widetilde G$ onto $V_N$.
Then the coefficients are uniquely determined by the orthogonality condition
\[
\int_0^1\int_0^1
\bigl(\widetilde G - P_t^{(N)} G\bigr)\cdot \psi
\,ds\,d\xi=0
\qquad (\text{for all}\ \psi\in V_N).
\]
Next, for a fixed integer $K\geq 0$ determined in the Bayesian optimization phase, we introduce the spatially band-limited space
\begin{equation*}
W_K:=
\biggl\{
\psi(t,s)=\sum_{|k|\le K}\psi_k(t)e^{2\pi i k s}
\biggr\}.
\end{equation*}
For each fixed $t$, we project in the $s$-direction onto the Fourier band $|k|\le K$:
\begin{equation*}
P_s^{(K)}f(t,s)
:=
\sum_{|k|\le K}\widehat f_k(t)e^{2\pi i k s},
\end{equation*}
where $\widehat f_k(t)$ denotes the $k$-th Fourier coefficient in $s$.
This restriction is natural when the design space itself is parameterized by finitely many Fourier modes, as in the Bayesian optimization phase.
Hence, the spatial filtering should be interpreted not merely as numerical smoothing, but as the orthogonal projection onto the prescribed design space.
Therefore the projected gradient is given by
\begin{equation*}
P_s^{(K)}P_t^{(N)}
=
P_t^{(N)}P_s^{(K)}:L^2\to V_N\cap W_K.
\end{equation*}
Hence, the practical update direction is
\begin{equation*}
\delta\Phi
=
\frac{P_s^{(K)}P_t^{(N)}G}
{\|P_s^{(K)}P_t^{(N)}G\|_{L^2}}.
\end{equation*}

\section{Numerical results}

In this section, we present numerical experiments to investigate
how the alignment between the vortex axis and the axis of swirling particles
is maintained during the time evolution of a vortex ring accompanied by the Kelvin wave.
In particular, we demonstrate that such alignment cannot be preserved by pure radial expansion alone, and that a nontrivial deformation
of the ring plays an essential role to keep alignment.
In the numerical implementation, we set
\begin{equation*}
\Gamma(t)=1-\cos(12\pi t),
\qquad
(t_0=1/48,\; t_1=1/24),
\end{equation*}
and adopt the following hyperparameters:
\begin{itemize}
\item $\epsilon = 0.2$ (strength of Kelvin wave)
\item $\omega = 240$ (temporal frequency of Kelvin wave)
\item $m =2$ (number of Kelvin wave)
\item $J = 1$ ($J+1$ polynomial terms, indexed by $j = 0, \dots, J$).

\end{itemize}
The time interval $[t_0,t_1]$ is discretized into $128$ steps,
and the angular parameter interval $[0,1)$ is discretized into $512$ points.
For the first step, we employ Bayesian optimization using the Optuna library (see \cite{Optuna}):
\begin{itemize}
\item First $500$ trials: quasi--Monte Carlo sampling
\item Subsequent $1{,}500$ trials: Bayesian optimization
\item Parameter bounds of $\{c^{\ell m}_{jk}\}$:
mean squared displacement at the terminal time $t_1$, more precisely,
$\sigma^2:=\sum_{k,\ell,m}|\sum_{j=0}^Jc_{jk}^{\ell m}|^2$  to be $\sigma = 0.3$
\end{itemize}

After the Bayesian optimization, we apply the
adjoint method:
\begin{itemize}
    \item $500$ iterations
    \item The learning rate is determined by an Armijo line search
    \item $N=8$ (projection for time direction)
\end{itemize}

\begin{figure}[!htbp]
\begin{center}
  \includegraphics[width=0.5\textwidth]{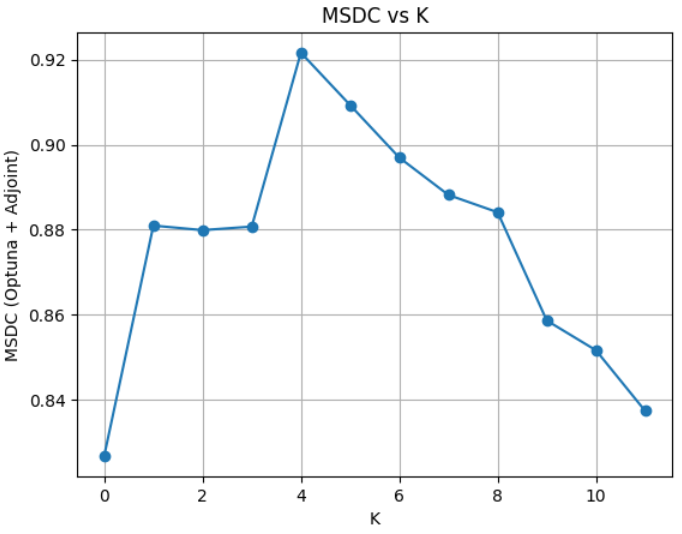}
\end{center}
\caption{MSDC}
\label{fig:MSDC}
\end{figure}

To demonstrate the scale locality, 
we vary $K$ ($K+1$ Fourier modes, indexed by $k=0,\cdots, K$), $K=0,1,\dots ,11$.
The numerical results are shown in Figure~\ref{fig:MSDC}.
We can observe the bump structure qualitatively,
and its peak is $K=4$.

\begin{figure}[!htbp]
\begin{center}
  \includegraphics[width=1.0\textwidth]{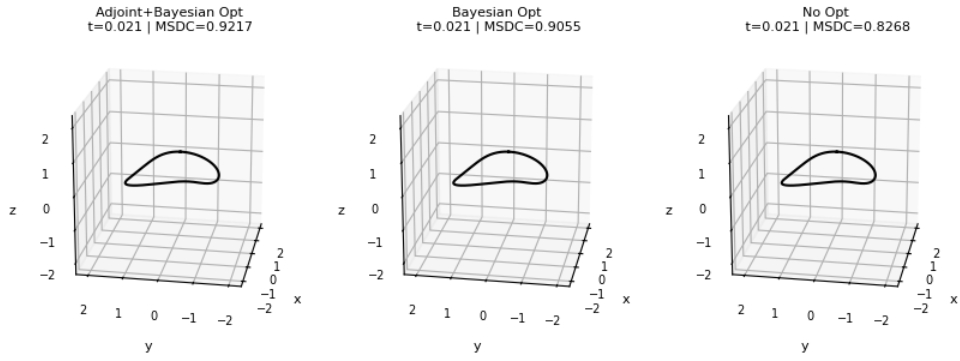}
\end{center}
\caption{Initial vortex rings (the same configuration). Left: Adjoint $+$ Bayesian optimizations, Center: Bayesian optimization, Right: No optimization.
}
\label{fig:initial}
\end{figure}
\begin{figure}[!htbp]
\begin{center}
  \includegraphics[width=1.0\textwidth]{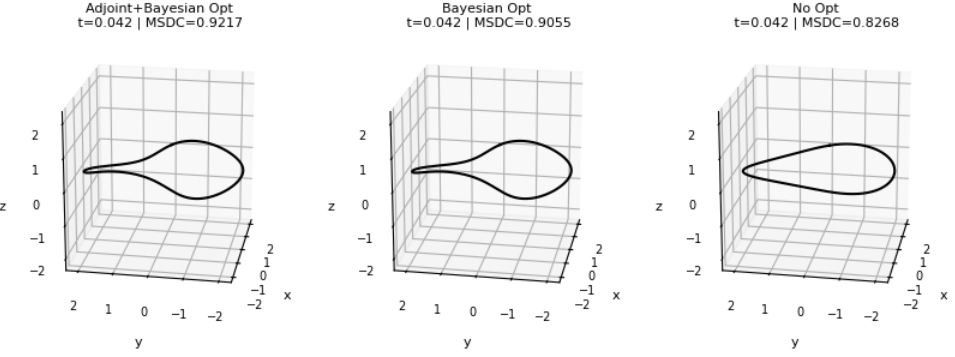}
\end{center}
\caption{Vortex rings at the terminal time. Left: Adjoint$+$ Bayesian optimizations, Center: Bayesian optimization, Right: No optimization.
}
\label{fig:terminal}
\end{figure}
For $K=4$, we visualize the results in Figures \ref{fig:initial} and~\ref{fig:terminal}.
The numerical experiments reveal that the present optimization problem exhibits a clear separation between global and local geometric components. 
A key observation is that adjoint-based optimization produces only marginal changes in the shape, 
while yielding a noticeable improvement in the objective functional (MSDC). 
From an optimization perspective, this implies that the objective landscape is highly non-convex and multi-modal with respect to the global geometric shape, while being relatively smooth within each basin.
As a consequence, adjoint-based methods are inherently local and can only explore a single basin of attraction, leading to invisible refinements without substantial shape deformation.
In contrast, global search methods (such as Bayesian optimization or quasi-Monte Carlo sampling) are capable of traversing multiple basins and identifying qualitatively different shapes.
These observations naturally justify a hybrid optimization strategy combining global exploration and local refinement:
\[
\text{global search} \;\longrightarrow\; \text{adjoint refinement}.
\]

Such combination of global and local optimization methods have been extensively studied in different contexts such as aerodynamic optimization (see \cite{Forrester2008, Jameson1988,Vicini1999} for example). On the other hand, to the best of our knowledge, such a combination of Lagrangian flow optimization and  differential geometry of curves
has received relatively limited attention.

\section{Conclusion}

Qualitatively, the numerical results indicate that maintaining alignment between the vortex axis and the axis of swirling particles throughout the evolution requires the presence of a nontrivial wavy deformation, and that the deformation frequency must not be too high, reflecting scale locality.  
At the present stage, our analysis remains qualitative, nevertheless, this observation provides strong evidence that the deformation of vortices (not only vortex rings) is governed by scale-local mechanisms.  
A central direction for future work is to undertake a detailed quantitative investigation of this scale locality for various vortex configurations.

\vspace{0.5cm}
\noindent
{\bf Acknowledgments.}\ 
The author thanks to Professors Susumu Goto, Chun Liu and Yoshikazu Giga for valuable comments.
Research of  TY  was partly supported by the JSPS Grants-in-Aid for Scientific
Research 24H00186.

\bibliographystyle{amsplain} 

\end{document}